\documentclass[12pt,a4paper,english]{amsart}
\usepackage{babel}
\usepackage[latin1]{inputenc}
\usepackage{amsmath}
\usepackage{amssymb}

\theoremstyle{plain}
\newtheorem{Thm}{Theorem}[section]

\newtheorem{Lemma}[Thm]{Lemma}

\theoremstyle{definition}
\newtheorem{Remark}[Thm]{Remark}

\newtheorem{Def*}{Definition}

\newcommand{\al}{\alpha}
\newcommand{\be}{\beta}
\newcommand{\ga}{\gamma}
\newcommand{\m}[1]{\mathrm{E}[#1]}

\begin{document}

\title[Completing a $k-1$ assignment]
{Completing a $k-1$ assignment}
\author[Svante Linusson]{Svante Linusson$^1$} 
\footnotetext[1]{Partially supported by EC's IHRP program through grant HPRN-CT-2001-00272}
\address{Svante Linusson, Dept of Mathematics\\ Link\"oping University\\ 581 83
Link\"oping\\ Sweden} \email{linusson@mai.liu.se}
\author[Johan W\"astlund]{Johan W\"astlund$^1$}
\address{Johan W\"astlund, Dept of Mathematics\\ Link\"oping University\\ 
581 83 Link\"oping\\ Sweden}
\email{jowas@mai.liu.se}


\date{\today}

\begin{abstract}
We consider the distribution of the value of the optimal
$k$-assignment in an $m \times n$-matrix, where the entries are independent 
exponential random variables with arbitrary rates.

We give closed formulas for both the Laplace transform of this random variable
and for its expected value  under the
condition that there is a zero-cost $k-1$-assignment.

\end{abstract}
\maketitle

\section{Introduction}

Let $M$ be an $m\times n$ matrix.
A $k$-{\em assignment} is a collection of 
$k$ entries in the matrix such that no two are in the same row or the 
same column. The {\em value} of a $k$-assignment is the sum of its entries. A $k$-assignment is called {\em optimal} if its value is no larger than the value of any other $k$-assignment. 
If the entries of the matrix $M$ are random variables then so is the value 
of the optimal $k$-assignment, here denoted ${\min}_{k}(M)$.

The study of the optimal $k$-assignment has been pursued by 
researchers from different fields and with different random variables 
as entries in $M$. The main focus has been to estimate the size of 
the expected value of ${\min}_{k}(M)$. For references and more details on 
the history see \cite{CS98} or \cite{LW03}. 

In 1998 Giorgio Parisi \cite{P98} conjectured that if $M$ is a 
$k\times k$ matrix 
with independent exponential random variables with rate 1, exp(1), 
then the expected value of the optimal $k$-assignment is
\[\m{{\min}_k(M)}=1+\frac{1}{4}+\frac{1}{9}+\dots +\frac{1}{k^2}.
\]

Two very different proofs of this conjecture were announced simultaneously 
in March 2003, \cite{LW03}\cite{NPS03}.
The beautiful conjecture of Parisi inspired much work on exact formulas and 
many different generalizations where studied
\cite{AS02,BCR02,CS98,CS02,EES01,LW00,N02}. In \cite{LW03} a formula 
for the expected value is given when the matrix entries are $\exp(1)$ or 0.

In this note we investigate the problem from a different extreme. We allow the rates of the exponential random variables to be arbitrary positive numbers. 
We include the infinity as a possible rate, which corresponds to
the entry being constant zero. 
We prove exact formulas for $\m{{\min}_k(M)}$ 
and for the Laplace transform $L({\min}_k(M))$ under the assumption that $\m{{\min}_{k-1}(M)}=0$.
Formulas for completing a $k-1$-assignment were considered previously 
in \cite{CS02} in the case when all rates are equal to 1. 
The Laplace transform for some special cases when the rates are 
all equal to 1, have been determined in \cite{AS02}.

We will first prove the slightly easier result on expected value and then 
compute the Laplace transform for the entire distribution. In theory one 
should be able to deduce the first result from the second but it seems 
easier to compute them separately.

\section{Preliminaries}
As in \cite{LW00,LW03} the concept of row and column covers of 
zeros will be important. The set of zeros will be denoted $Z$.
We will consider sets of rows and columns in the $m\times
n$-matrix. A set $\al$ of rows and columns is said to {\em cover}
$Z$ if every entry in $Z$ is on either a row or on a column in
$\al$. A cover with $k-1$ rows and columns will be called a
$k-1$-cover. For many readers it might be convenient to translate the matrix 
to a bipartite graph. In that setting our covers are so called vertex covers.

Given a set of rows and columns $\al$ let the {\em rectangle} $R(\al)$ be the part of the matrix not covered by $\al$. If $\al$ is a $k-1$-cover of the zeros then the corresponding rectangle $R(\al)$ will be called {\em critical}. 
Let $Q_{k,Z}$ be the set of all critical rectangles in $M$.
We define a partial ordering on $Q_{k,Z}$ by letting,
for $R(\al),R(\be)\in Q_{k,Z}$, $R(\al) \leq R(\be)$ if 
the set of columns in $\al$ is a subset of 
the set of columns in $\be$, and the set of rows in $\al$ is a superset of the 
set of rows in $\be$. 

Recall that a random variable $X \sim \exp(\alpha)$ if $P(X>t)=e^{-t\alpha}$. Then $\alpha$ is called the {\em rate} and $\m{X}=1/\alpha$. 
If $X_1,\dots, X_n$ are independent and $X_i\sim \exp(\alpha_i)$ then $\m{\min_i X_i}=1/\sum_i \alpha_i$.
Let the {\em rate} of a rectangle $R$ be the
sum of the rates of the individual entries.
We denote it $I(R)$ and note that $0<I(R)\le \infty$. 

From matching theory we have the following lemma that 
will be important.
\begin{Lemma}
Assume that $Z$ contains $k-1$ zeros in independent position. Then 
$Q_{k,Z}$ is a lattice. In particular it has unique maximal and 
minimal elements.
\end{Lemma}

{}From Theorem 2.9 in \cite{LW00} we cite the following useful fact.

\begin{Lemma}\label{L:row}
Suppose that a row $r$ belongs to every $k-1$-cover of the zeros $Z$ in the 
matrix $M$. Then every optimal $k$-assignment contains a zero from row $r$. 
This means that we can remove row $r$ from $M$ and obtain a matrix 
$M^r$ with the 
property that $\m{{\min}_{k}(M)}= \m{{\min}_{k-1}(M^r)}$. 
\end{Lemma}

Recall that the {\em incidence algebra} over a poset is the algebra of 
functions defined on the intervals in the poset, see \cite{S}. The product 
in the incidence algebra is defined by convolution. One way to visualize 
the incidence algebra over a poset $Q$ is to let a function $f$ be 
represented by a matrix $F$ with rows and columns indexed by the elements in $Q$. The element on position $(\al,\be)$ 
in $F$ is $f(\al,\be)$. The rows and columns must be arranged according to a linear extension of $Q$. Every such matrix is upper triangular and multiplication of functions corresponds to multiplication of matrices.

\section{The expected value}

Define a function in the incidence algebra of $Q_{k,Z}$ by 
$f(R(\al),R(\be))=I(R(\al)\cap R(\be))$. Since $I(R(\al))> 0$ for 
all critical rectangles, $f$ is invertible in the incidence algebra.

\begin{Thm}\label{T:expected}
Let $M$ be an $m\times n$ matrix with all entries being either zero or independent exponential random variables with arbitrary positive rates. 
Assume that the set of zeros
$Z$ in $M$ contains $k-1$ zeros in independent position, i.e.
$\m{{\min}_{k-1}(M)}=0$.

Then

\[ \m{{\min}_{k}(M)}=\sum_{R(\al)\le R(\be)} f^{-1}(R(\al),R(\be)),
\]
where the sum is over all intervals in $Q_{k,Z}$ and $f^{-1}$ is the inverse of 
$f$ in the incidence algebra. Equivalently we can write

\begin{align}
&\m{{\min}_{k}(M)}=\notag\\
&\sum_{R_1<R_2<\dots<R_s}(-1)^{(s-1)} \frac{I(R_1 \cap R_2)
I(R_2 \cap R_3)\cdots I(R_{s-1} \cap R_s)}{I(R_1)I(R_2)\cdots I(R_s)},\notag
\end{align}
where the sum is taken over all non-empty chains in $Q_{k,Z}$. 
\end{Thm}

\begin{Remark} \label{R:comp}
The second formula runs over all chains in $Q_{k,Z}$ which 
in the worst case has size of order $k!$. 
The first formula is a computational and conceptual 
improvement for large $k$. 
It involves taking the inverse of a matrix indexed by the 
elements of $Q_{k,Z}$ which is exponential in $k$ in the worst case.

\end{Remark}

\begin{proof} The equivalence follows from Lemma \ref{L:inverse} below.
The proof will be by induction over $k$. The theorem is certainly true
for $k=1$.
Without loss of generality we may assume that 
the entries $(1,1),(2,2),\dots (k-1,k-1)$ are zero entries.
We may also assume that 
the maximal rectangle $R(\ga)$ in $Q_{k,Z}$ corresponds to 
the cover $\ga$ consisting of columns $1,\dots,k-1$.
If this is not the case, then there is a row $i$ that belongs to every 
$k-1$-cover. This implies by Lemma \ref{L:row} that 
$\m{{\min}_{k}(M)}= \m{{\min}_{k-1}(M^i)}$, where $M^i$ is obtained from $M$ 
by completely removing row $i$. Since these matrices have the same $Q_{k,Z}$ 
the result is clear by induction.

Note that any chain in $Q_{k,Z}$ not ending with $R(\ga)$
can be augmented with $R(\ga)$ at the end.

We now use the same recursion procedure as in \cite{AS02} and \cite{LW00}
corresponding to the optimal cover $\ga$ consisting of the first $k-1$ 
columns.
That is, let $X$ be the minimum of all the exponential random variables 
in $R(\ga)$. Note that exactly one entry 
in $R(\ga)$ will belong to an optimal $k$-assignment and that $\m{X}=1/I(R(\ga))$. Subtract $X$ from all entries in $R(\ga)$ and a 
new zero will occur at the position of the minimum. All other entries 
will be unchanged in distribution by the forgetfulness of the exponential distribution.
Let $M_{i,j}$ be the matrix obtained when position $(i,j)$ in $M$ has been replaced with a zero. 

Let $K_i$ be the intersection of $R(\ga)$ and row $i$.
If the minimum is in $K_i$, $i>k-1$, then we get a zero $k$-assignment in $M_{i,j}$.
If the minimum is in $K_i$, $i=1,\dots,k-1$,
then row $i$ has to be in every optimal cover of $M_{i,j}$ and as above 
we can remove row $i$ and this case is done by induction. 
Let again $M^i$ denote the matrix with row
$i$ removed from $M$. Also let $R^i$ be the column maximal cover of $M^i$.
The poset of $k-1$-covers of $M^i$ will be the induced subposet of 
$Q_{k,Z}$ of elements $\le R^i$.

The probability that the zero occurs in $K_i$ is 
$\frac{I(K_i)}{I(R(\ga))}$ and we get
\[\m{{\min}_{k}(M)}=\frac{1}{I(R(\ga))}+\sum_{i=1}^{k-1}\frac{I(K_i)}{I(R(\ga))} 
\m{{\min}_{k-1}(M^i)}.
\]

Which by induction becomes

\[\frac{1}{I(R(\ga))}+\sum_{i=1}^{k-1}\frac{I(K_i)}{I(R(\ga))}
\sum_{R_1<\dots<R_s\le R^i}(-1)^{(s-1)} \frac{I(R_1 \cap R_2)
\cdots I(R_{s-1} \cap R_s)}{I(R_1)\cdots I(R_s)}.
\]

Change the order of summation to get
\begin{align}
&\frac{1}{I(R\ga)}+\notag\\
&\sum_{R_1<\dots<R_s<R(\ga)}(-1)^{(s-1)} \frac{I(R_1 \cap R_2)
\cdots I(R_{s-1} \cap R_s)}{I(R_1)\cdots I(R_s)}
\sum_{i\notin rowset(R_s)}\frac{I(K_i)}{I(R(\ga))}.\notag
\end{align}
Here $rowset(R)$ denotes the set of rows of $M$ 
that intersect the rectangle $R$.
Now 
\[
\sum_{i\notin rowset(R_s)}\frac{I(K_i)}{I(R(\ga))}=
\frac{I(R(\ga)\backslash R_s)}{I(R(\ga))}=1-
\frac{I(R_s\cap R(\ga))}{I(R(\ga))},
\]
since the entries are independent exponential random variables.
The theorem follows.
\end{proof}

\section{The Laplace transform}

We can in fact use the same proof technique to get the stronger result 
determining the Laplace transform of the distribution of ${\min}_{k}(M)$.

Recall that the Laplace transform of a random variable $X$ is 
$L(X,t)=\m{e^{-tX}}$. It has the following well-known properties.

\[L(X+Y,t)=L(X,t)L(Y,t), \qquad\text{if $X$ and $Y$ are independent}\]

Given random variables $X_1,\dots, X_s$ and probabilities 
$p_1,\dots, p_s$, define the random variable $I$ to take value $i$ 
with probability $p_i$, independent of $X_1,\dots,X_s$. Then 
\[L(X_I,t)=\sum_{i=1}^s p_iL(X_i,t).
\]
In this situation we will need the special case $L(0,t)=1$.

For a critical rectangle $R$ we will use the notation 
\[
\phi(R,t)=L({\min}_{1}(R))=\frac{I(R)}{I(R)+t}.\]

As for Theorem \ref{T:expected} we give two statements of the same formula
using Lemma \ref{L:inverse}. Remark \ref{R:comp} applies also here. 

\begin{Thm}
Let $M$ be an $m\times n$ matrix with all entries being either zero or 
exponential independent random variables with arbitrary positive rates. 
Assume that the set of zeros
$Z$ in $M$ contains $k-1$ zeros in independent position, i.e.
$\m{{\min}_{k-1}(M)}=0$.

Then we can write 
$L({\min}_{k}(M))=\sum_{R\in Q_{k,Z}} c_R(M)\phi(R,t)$. Furthermore 

$c_R(M)=a_R\cdot b_R,$ where

\[a_R=\sum_{R_1<R_2<\dots<R_s=R}(-1)^{s} \prod_{i=1}^{s-1}
\frac{I(R_i \cap R_{i+1})-I(R)}{I(R_i)-I(R)}\] and
\[b_R=\sum_{R=R_s<R_{s+1}<\dots<R_u}(-1)^{(u-s+1)} \prod_{i=s+1}^{u}
\frac{I(R_i \cap R_{i-1})-I(R)}{I(R_i)-I(R)},
\]
where the sums are taken over all chains containing $R$ in $Q_{k,Z}$.

In the generic case rewrite this as
\[
c_R(M)=\left(\sum_{R(\al)\le R} g_R^{-1}(R(\al),R)\right)\cdot
\left(\sum_{R\le R(\be)} g_R^{-1}(R,R(\be))\right),
\]
where the sums are over elements in $Q_{k,Z}$ and where $g^{-1}$ 
is the inverse in the incidence algebra of 
\[g_R(R_i,R_j)=\begin{cases} 
1, &\text{if $R_i=R_j=R$}\\
0, &\text{if $R_i\nleq R_j$}\\
I(R_i\cap R_j)-I(R), &{\text otherwise.}
\end{cases}
\] 
\end{Thm}

\begin{proof} The proof is by induction over $k$. We use the same 
notations as in the proof of Theorem \ref{T:expected} and compute 
the Laplace transform using the same recursive step.

\[ L({\min}_k(M),t)=\phi(R(\ga),t)\left( p+\sum_{i=1}^{k-1}\frac{I(K_i)}{I(R(\ga))}
\cdot L({\min}_k(M^i)) \right),
\]
where $p$ is the probability that the minimum is located so a zero cost $k$-assignment occurs. We now decompose this product by the method 
of partial fractions with respect to $t$. This proves the first claim and the 
uniqueness of the coefficients $c_R(M)$. The function $g_R$ is invertible if
$I(R)\neq I(R(\al))$ for all $R(\al)\in Q_{k,Z}$ such that $R\le R(\al)$ or
$R(\al)\le R$. The equivalence of the 
two formulas for the coefficients then follows from Lemma \ref{L:inverse}.  

After the decomposition by partial fractions 
follows an extraction of the terms involving $\phi(R,t)$ which gives
\[c_R(M)=\frac{I(R(\ga))}{I(R(\ga))-I(R)}\sum_{i=1}^{k-1}
\frac{I(K_i)}{I(R(\ga))}\cdot c_R(M^i)
\]
Assume that $R\neq R(\ga)$ and we may inductively write
\begin{align}
c_R(M)=&\sum_{i=1}^{k-1}
\frac{I(K_i)}{I(R(\ga))-I(R)} a_R(M) \cdot\notag\\
&\sum_{R=R_s<R_{s+1}<\dots<R_u}(-1)^{(u-s+1)} \prod_{i=s+1}^{u}
\frac{I(R_i \cap R_{i-1})-I(R)}{I(R_i)-I(R)},\notag
\end{align}
where $R_u$ does not intersect the row $i$. Change the order of summation to get
\begin{align}
c_R(M)= a_R(M)
&\sum_{R=R_s<R_{s+1}<\dots<R_u}(-1)^{(u-s+1)} \cdot\notag\\
&\prod_{i=s+1}^{u} \frac{I(R_i \cap R_{i-1})-I(R)}{I(R_i)-I(R)}\cdot \sum_{i\notin colset(R_u)}
\frac{I(K_i)}{I(R(\ga))-I(R)},\notag
\end{align}
and proceed as in the proof of Theorem \ref{T:expected}. 
Here $colset(R)$ denotes the set of columns of $M$ 
that intersect the rectangle $R$. Note that trivially
$a_R(M)=a_R(M^i)$. If $R=R(\ga)$ the recursive step uses 
the row maximal rectangle instead, that is the 
smallest element in $Q_{k,Z}$ to compute $a_R(M)$.
\end{proof}

\section{A lemma}

\begin{Lemma} \label{L:inverse}
Let $f$ be a function in the incidence algebra over a poset $P$ and let 
$\al \le \be\in P$ be arbitrary elements. The inverse of $f$ can be written as
\[ f^{-1}(\al,\be)=
\sum_{\al=\ga_1<\ga_2<\dots<\ga_s=\be}(-1)^{(s-1)} \frac{f(\ga_1,\ga_2)
f(\ga_2,\ga_3)\cdots f(\ga_{s-1},\ga_s)}{f(\ga_1,\ga_1)f(\ga_2,\ga_2)\cdots 
f(\ga_s,\ga_s)},
\]
where the sum is over all chains in the interval $[\al,\be]$ 
beginning in $\al$ and ending in $\be$.
\end{Lemma}

\begin{proof} Let $F$ be the upper triangular matrix corresponding to $f$ as 
described in the preliminaries. 
Let $D$ be the diagonal matrix with the values of 
$f(\ga,\ga)$ on the diagonal and zeroes elsewhere. Let $N$ be the nilpotent 
matrix that agrees with $F$ at all positions except on the diagonal, 
where $N$ has zeros. We can then write $F=D+N$ and we can easily verify 
that 
\[F^{-1}=D^{-1}- D^{-1}ND^{-1}+ D^{-1}ND^{-1}ND^{-1} - \dots.
\]
The matrix $ D^{-1}N $ will have zeros on and below the diagonal. Thus the sum is finite. Since inverting $f$ in the incidence algebra is the 
same as inverting the matrix $F$ the lemma follows.
\end{proof}

\begin{Remark}
The main theorem in \cite {LW03} when the exponential random variables all
have rate 1 or infinity has a reformulation in  
\cite{LW00} in terms of the M\"obius function of a certain poset 
called $P$. 
That poset is different and 
much larger than $Q_{k,Z}$. To be more precise the atoms in $P$ 
are the elements in $Q_{k,Z}$. All our efforts to join the two theorems to one
for completely arbitrary rates have so far been fruitless.
\end{Remark}

\end{document}